\documentclass[12pt,reqno]{amsart}
	\usepackage{amssymb}
	\usepackage[mathscr]{eucal}
	\usepackage{xypic}
	\usepackage{amsmath}
	\usepackage{epsfig}
	\usepackage{amscd}
	\usepackage{mathrsfs}
\usepackage[english]{babel}

	\usepackage{tikz-cd}
	\usepackage{hyperref}
	\theoremstyle{definition}
	\newtheorem{theorem}{Theorem}[section]
	\newtheorem{definition}{Definition}[section]
	\newtheorem{lemma}[theorem]{Lemma}
	\newtheorem{proposition}[theorem]{Proposition}
	\newtheorem{corollary}[theorem]{Corollary}
	\newtheorem{remark}[theorem]{Remark}
	\newtheorem{example}[theorem]{Example}

	\newcommand{\mc}{\mathcal}
	\newcommand{\mb}{\mathbb}
	\newcommand{\mf}{\mathfrak}
	\newcommand{\xra}{\xrightarrow}
	\newcommand{\ra}{\rightarrow}
	\newcommand{\rra}{\rightrightarrows}
	\numberwithin{equation}{section}

	\def\og{\leavevmode\raise.3ex\hbox{$\scriptscriptstyle\langle\!\langle$~}}
	\def\fg{\leavevmode\raise.3ex\hbox{~$\!\scriptscriptstyle\,\rangle\!\rangle$}}
	
\begin{document}
		
		\title[Extension of topological groupoids and Serre, Hurewicz morphisms]
		{Extension of topological groupoids and Serre, Hurewicz morphisms}
		
		\author[S. Chatterjee]{Saikat Chatterjee}
		
		\address{School of Mathematics,
			Indian Institute of Science Education and Research--Thiruvananthapuram,
			Maruthamala P.O., Vihtura, Kerala 695551, India}
		\email{saikat.chat01@gmail.com}
		
		\author[P. Koushik]{Praphulla Koushik}
		
		\address{School of Mathematics,
			Indian Institute of Science Education and Research--Thiruvananthapuram,
			Maruthamala P.O., Vihtura, Kerala 695551, India}
		\email{koushik16@iisertvm.ac.in}

		\subjclass[2010]{Primary 18F20, Secondary 22A22, 53C08}
		
		\keywords{Topological groupoid extensions; Gerbes; Serre, Hurewicz morphisms}
	\begin{abstract}
		In this paper, we introduce the notion of a topological groupoid extension and relate it to the already existing notion of a gerbe over a topological stack. We further study the properties of a gerbe over a Serre, Hurewicz stack.
	\end{abstract}
	\maketitle
	
	\section{Introduction}
	Related areas of topological groupoids  are quite an active area of research. Here we refer to only a few of them, such as  \cite{MR2927363}, \cite{MR2977576}, \cite{MR1633759}, \cite{MR1041882}. However, the topological stacks and related geometry are not explored to a great extent, other than  by Behrang Noohi in	\cite{noohi2005foundations}, \cite{MR3144243},  \cite{MR2719557}, \cite{MR2977576}, \cite{MR2927363}, \cite{MR3552548} and some other authors  \cite{MR2900442},	\cite{ebert2009homotopy}, \cite{metzler2003}. On the other hand, in the smooth set-up the Morita equivalence of Lie groupoid extensions and their association to differentiable stacks and related geometry have been  studied in several articles, for example  \cite {MR2493616}, \cite{MR4124773}, \cite{MR2778793}. 	
	A version of Morita equivalence of groupoids in the topological set-up has been introduced in the thesis of Carchedi \cite{Carchedi}. In this paper, we are going to introduce the notion of Morita equivalence of topological groupoid extensions and explore its relation with a gerbe over a topological stack. Our main objective is to study the properties of such a gerbe over Serre, Hurewicz stacks. The notion of Serre, Hurewicz stacks was introduced in \cite{MR3144243} by Noohi, and this paper is the main motivation behind this work. 
	
	Among the other papers referenced above, \cite{noohi2005foundations} and \cite{MR2927363} discussed the homotopy theories of topological stacks. 	In \cite{MR2977576}, the authors explored the relation between topological stacks and string topology. In \cite{MR2719557} Noohi studied the set of maps between topological stacks. In particular the paper shows  $\text{Hom}(\mc{X},\mc{Y})$ is a topological stack if $\mc{Y}$ has a groupoid representation by a compact topological space. Whereas	\cite{ebert2009homotopy} studied the topological stacks representable by paracompact groupoids. The definition of a gerbe over a topological stack was introduced in \cite{metzler2003}, which we are going to adopt in this article. For the discussions on the notions of topological groupoid actions and principal bundles in the set-up of sites, we refer to \cite{MR3438234}.	 Some of the other papers which we have consulted are \cite{MR2292996}, \cite{MR2408155},  \cite{MR973173}, \cite{MR413096}, \cite{MR2729658}.

	Let ${\rm Top}$ be the category of topological spaces. We fix the open cover Grothendieck topology on ${\rm Top}$.
	We define a topological groupoid $\mc{G}=[\mc{G}_1\rra \mc{G}_0]$ to be a groupoid internal to the category of topological spaces ${\rm Top}$. 
	A $\mc{G}$-principal bundle or a $\mc{G}$-torsor  is a map of local sections $\pi\colon P\,\ra M$ with an action of $\mc{G}$ on $P$ satisfying certain conditions. Then it is a common knowledge  that the category of $\mc{G}$ bundles $B\mc{G}$ defines a stack over ${\rm Top}$. In particular a stack over ${\rm Top}$ is called a topological stack if it is isomorphic to $B\mc{G}$ for some topological groupoid $\mc{G}$. A gerbe is basically a morphism of topological stacks $F\colon \mc{C}\ra\mc{D}$ with some additional properties (namely both $F\colon  \mc{D}\ra \mc{C}$ and the diagonal  morphism $\Delta_F\colon  \mc{D}\ra \mc{D}\times_{\mc{C}}\mc{D}$ are epimorphism of topological stacks)~\cite{metzler2003}. A topological groupoid extension is a morphism of topological groupoids $(F,1_M)\colon  [\mc{G}_1\rra M]\ra [\mc{H}_1\rra M]$ such that $F\colon  \mc{G}_1\ra \mc{H}_1$ admits local sections. 
	We call a morphism of topological groupoids $[\Delta_1\rra \Delta_0]\ra [\Gamma_1\rra \Gamma_0]$  Morita morphism of topological groupoids if $\Delta_0\ra \Gamma_0$ is a map of local sections and $[\Delta_1\rra \Delta_0]$ is the pullback groupoid of $[\Gamma_1\rra \Gamma_0]$ along the morphism $\Delta_0\ra \Gamma_0$. If there exists a pair of Morita morphisms $[P_1\rra P_0]\ra [\Delta_1\rra \Delta_0]$ and $[P_1\rra P_0]\ra [\Gamma_1\rra \Gamma_0]$, then we say $[\Delta_1\rra \Delta_0]$ and $[\Gamma_1\rra \Gamma_0]$ are Morita equivalent.  This definition of Morita equivalence naturally extends to the definition of  Morita equivalence of topological groupoid extensions. 
	
	In this paper, we show that a Morita equivalence class of topological groupoid extension defines a gerbe over a topological stack. On the other hand a gerbe $F\colon \mc{C}\ra\mc{D}$ for which the diagonal  morphism $\Delta_F\colon  \mc{D}\ra \mc{D}\times_{\mc{C}}\mc{D}$ is a representable map of local sections defines a Morita equivalence class of topological groupoid extensions. It should be noted here that a representable map of local sections is a slightly stronger condition than an epimorphism. The corresponding result in smooth set-up has been proven by the authors in \cite{MR4124773}.

	
	Below we give a brief outline of this article.
	
	Section \ref{Section:top groupoids and top stacks} and \ref{Section:morphism of top groupoids and extensions} 
	are expository in nature. Though we introduce some new definitions, such as extensions and  Morita equivalences of topological groupoids in Section~\ref{Section:top groupoids and top stacks}, these sections mainly review and recall the already existing definitions and results. In Section~\ref{Section:morphism of top groupoids and extensions}	we outline the construction of a morphism of stacks associated to a map of topological groupoids and recall the notion of a gerbe.   In Section \ref{section:Topogrpdextensionmorstacks} we prove the correspondence between topological groupoid extensions and gerbes. In Section~\ref{section:Serre-hurewiczfibrations} we show that if $[\mc{G}_1\rra M]$ is a topological groupoid with  the source map $s\colon  \mc{G}_1\ra M$  a locally Hurewicz fibration (respectively, Serre fibration), then the gerbe associated to an extension $[\mc{G}_1\rra M]\,\ra \, [\mc{H}_1\rra M]$ is a Hurewicz gerbe (respectively, a Serre gerbe).
	
	\section{Topological groupoids and topological stacks}\label{Section:top groupoids and top stacks}
	We begin with the simplest and most general definition of a topological groupoid. We will mostly rely on the definitions in { \cite{Carchedi}}.
	By a \textit{topological groupoid}, we mean a groupoid $\mc{G}$ whose object and morphism sets are topological spaces, and all structure maps are continuous.  We denote the object and morphism sets, respectively, by $\mc{G}_0$ and $\mc{G}_1$. Some other authors have defined a topological groupoid requiring the structure maps to satisfy some additional conditions and studied properties of such topological groupoids. It would also be our interest to study some of those properties. However, instead of incorporating the conditions in the definition itself, we will impose the conditions when and where required. \textit{A morphism $\phi\colon\mc{G}\,\ra\, \mc{H}$ between topological groupoids}  is a continuous functor compatible with all the structure maps of the groupoids.  	We say a continuous map $f\colon X\ra Y$ \textit{admits local sections} or is  a \textit{map of local sections} if there is an open cover $\{U_{\alpha}\}$ of $Y$ and a continuous section $U_{\alpha}\,\ra\, X$ for each $\alpha$.

	\begin{definition}[topological groupoid extension]	
		Let $[\mc{H}_1\rra M]$ be a topological groupoid. A \textit{topological groupoid extension of $[\mc{H}_1\rra M]$} is a morphism of topological groupoids $(F,1_M)\colon [\mc{G}_1\rra M]\ra [\mc{H}_1\rra M]$ such that $F\colon  \mc{G}_1\ra \mc{H}_1$ admits local sections. We often denote such a topological groupoid extension by $F\colon  \mc{G}_1\ra \mc{H}_1\rra M$.\end{definition}
	
	The Morita equivalence of Lie groupoids and Lie groupoid extensions have very interesting properties and have been studied in ~\cite{MR2493616} and \cite{MR4124773}. The counterpart of the same in the topological set-up has not been studied. In fact, even we are not aware of any existing definition of  Morita equivalent topological groupoid extensions. Here we introduce the notion of Morita equivalent topological groupoid extensions and study some of their properties. Particularly we are interested in its association with the topological gerbes.  
	
	Let $[\Gamma_1\rra \Gamma_0]$ be a topological groupoid and $J:P_0\ra \Gamma_0$ be a map of local sections. Let $P_1$ denote the fiber product 
	$P_0\times_{J,\Gamma_0,s}\Gamma_1\times_{t,\Gamma_0,J}P_0$. Then $[P_1\rra P_0]$ has a topological groupoid structure with source and target maps being the projection maps ${\rm pr}_1, {\rm pr}_3\colon P_1\ra P_0$ respectively. We call $[P_1\rra P_0]$ to be the \textit{pullback groupoid} of $[\Gamma_1\rra \Gamma_0]$ along $J\,\colon\,P_0\ra \Gamma_0$.
	
	A morphism of topological groupoids $[\Delta_1\rra \Delta_0]\ra [\Gamma_1\rra \Gamma_0]$ is said to be a \textit{Morita morphism of topological groupoids} if $\Delta_0\ra \Gamma_0$ is a map of local sections and $[\Delta_1\rra \Delta_0]$ is isomorphic to the pullback groupoid of $[\Gamma_1\rra \Gamma_0]$ along the morphism $\Delta_0\ra \Gamma_0$. 
	
	\begin{definition}
		Let $[\Gamma_1\rra \Gamma_0]$ and $[\Delta_1\rra \Delta_0]$ be a pair of  topological groupoids. We say that $[\Gamma_1\rra \Gamma_0]$ and $[\Delta_1\rra \Delta_0]$ are \textit{Morita equivalent} if there exists a third topological groupoid $[P_1\rra P_0]$ with a pair of Morita morphisms 
		$[P_1\rra P_0]\ra [\Delta_1\rra \Delta_0]$ and $[P_1\rra P_0]\ra [\Gamma_1\rra \Gamma_0]$.
	\end{definition}
	\begin{definition}[Morita morphism of topological groupoid extensions]\label{Definition:moritamorphismtopologicalextension}
		Let $\phi'\colon X_1'\rightarrow  Y_1'\rightrightarrows M'$  and 
		$\phi\colon X_1\rightarrow  Y_1\rightrightarrows M$ be a pair of  topological groupoid extensions.
		A \textit{Morita morphism of topological groupoid extensions} from   $\phi'\colon X_1'\rightarrow  Y_1'\rightrightarrows M'$  to 
		$\phi\colon X_1\rightarrow  Y_1\rightrightarrows M$ is given by a pair of Morita morphisms of topological groupoids, 
		\begin{equation}\begin{tikzcd}[sep=small]
		X_1' \arrow[dd,xshift=0.75ex,"t"]\arrow[dd,xshift=-0.75ex,"s"'] \arrow[rr, "\psi_X"] &  & X_1 \arrow[dd,xshift=0.75ex,"t"]\arrow[dd,xshift=-0.75ex,"s"'] \\
		&  &  \\
		M' \arrow[rr, "f"] &  & M
		\end{tikzcd}\text{ and }\begin{tikzcd}[sep=small]
		Y_1' \arrow[dd,xshift=0.75ex,"t"]\arrow[dd,xshift=-0.75ex,"s"'] \arrow[rr, "\psi_Y"] &  & Y_1 \arrow[dd,xshift=0.75ex,"t"]\arrow[dd,xshift=-0.75ex,"s"'] \\
		&  &  \\
		M' \arrow[rr, "f"] &  & M
		\end{tikzcd}\end{equation}
		such that the following diagram
		\begin{equation}\begin{tikzcd}[sep=small]
		X_1' \arrow[dd,"\psi_X"'] \arrow[rr, "\phi'"] &  & Y_1' 
		\arrow[dd,"\psi_Y"] \\
		&  &  \\
		X_1 \arrow[rr, "\phi"] &  & Y_1
		\end{tikzcd}\end{equation} is commutative.
	\end{definition}
	\begin{definition}[Morita equivalent topological groupoid extensions]\label{Definition:moritaequivalenttopologicalgroupoidextensions}
		Let $\phi'\colon X_1'\rightarrow  Y_1'\rightrightarrows M'$  and 
		$\phi\colon X_1\rightarrow  Y_1\rightrightarrows M$ be a pair of topological groupoid extensions. We say that  $\phi'\colon X_1'\rightarrow  Y_1'\rightrightarrows M'$  and 
		$\phi\colon X_1\rightarrow  Y_1\rightrightarrows M$ are  \textit{Morita equivalent topological groupoid extensions}, if there exists a third topological groupoid extension $\phi''\colon X_1''\rightarrow Y''\rightrightarrows M''$ and a pair of Morita morphisms of topological groupoid extensions
		\[(\phi''\colon X_1''\rightarrow Y''\rightrightarrows M'')\rightarrow (\phi\colon X_1\rightarrow Y_1\rightrightarrows M)\] and \[(\phi''\colon X_1''\rightarrow Y''\rightrightarrows M'')\rightarrow (\phi'\colon X_1'\rightarrow Y_1'\rightrightarrows M').\]
	\end{definition}
	
	\begin{definition}[{\cite
			{Carchedi}}]\label{Definition:ActionOfTopGroupoid}
		Let $\mc{G}=[\mc{G}_1\rra \mc{G}_0]$ be a topological groupoid. Let $P$ be a topological space. A \textit{left action of the topological groupoid $\mc{G}$ on the topological space $P$} is given by a pair of continuous maps $(a_{\mc{G}}\colon P\ra \mc{G}_0,\mu\colon \mc{G}_1\times_{s,\mc{G}_0,a_{\mc{G}}}P\ra P)$
		satisfying the following conditions:
		\begin{enumerate}
			\item $a_{\mc{G}}(\mu(\gamma,p))=t(\gamma)$ for all $(\gamma,p)\in \mc{G}_1\times_{s,\mc{G}_0,a_{\mc{G}}}P$,
			\item $\mu(1_{a_{\mc{G}}(p)},p)=p$ for all $p\in P$,
			\item $\mu(\gamma,\mu(\gamma',p))=
			\mu(\gamma\circ \gamma',p)$ for all $(\gamma,\gamma',p)\in \mc{G}_1\times_{s,\mc{G}_0,t}\mc{G}_1\times_{s,\mc{G}_0,a_{\mc{G}}}P$.
		\end{enumerate} 
		
		Similarly  a \textit{right action} is a pair 
		of continuous maps $(a_{\mc{G}}\colon P\ra \mc{G}_0,\mu\colon P\times_{a_{\mc{G}},\mc{G}_0,t}\mc{G}_1\ra P)$ satisfying 
		\begin{enumerate}		
			\item $a_{\mc{G}}(\mu(p,\gamma))=s(\gamma)$ for all $(p,\gamma)\in P\times_{a_{\mc{G}},\mc{G}_0,t}\mc{G}_1$,
			\item $\mu(p,1_{a_{\mc{G}}(p)})=p$ for all $p\in P$,
			\item $\mu(\mu(p,\gamma),\gamma')=\mu(p, \gamma\circ \gamma')$ for all $(p,\gamma,\gamma')\in P\times_{a_{\mc{G}},\mc{G}_0,t}\mc{G}_1\times_{s,\mc{G}_0,t}\mc{G}_1$.
		\end{enumerate} 	
	\end{definition}
	Note that the topological spaces $\mc{G}_1\times_{s,\mc{G}_0,a_{\mc{G}}}P$ and $P\times_{a_{\mc{G}},\mc{G}_0,t}\mc{G}_1$ will respectively be denoted as $\mc{G}_1\times_{\mc{G}_0}P$
	and $P\times_{\mc{G}_0}\mc{G}_1$. We will also write   $p\cdot \gamma$ instead of   $\mu(p, \gamma)$.
	
	\begin{definition}[{\cite{Carchedi}}]\label{Definition:G-bundle}
		Let $\mc{G}=[\mc{G}_1\rra \mc{G}_0]$ be a topological groupoid and $M$ be a topological space. 
		A \textit{principal $\mc{G}$-bundle over the topological space $M$} is a topological space $P$ with a   continuous map  $\pi\colon P\ra M$ admitting local sections and a right action $(a_{\mc{G}}, \mu)$ of 	$\mc{G}$ on $P$ such that 		
		\begin{enumerate}
			\item $\pi(p\cdot\gamma)=\pi(p)$ for all $(p,\gamma)\in P\times_{\mc{G}_0}\mc{G}_1$,
			\item the map $P\times_{\mc{G}_0}\mc{G}_1\ra P\times_{M}P$, $(p,\gamma)\mapsto (p,p\cdot \gamma)$ is a homeomorphism.
		\end{enumerate}
	\end{definition} 
	Observe that the condition (1) above ensures that the action is fiber preserving, whereas (2)  ensures that the action is transitive and free on the fibers. Note that in \cite{Carchedi}, the $\mathcal{G}$ acts on the left, whereas, here we have considered as a right action.
	\begin{definition}\label{Def:mapG-bundle}Let $\mc{G}=[\mc{G}_1\rra \mc{G}_0]$ be a topological groupoid and $(P,\pi,M), (P',\pi',M')$ principal $\mc{G}$-bundles. A \textit{morphism} from $P(\pi,M)$ to $P'(\pi',M')$ is given by a pair of continuous maps
		$(F\colon P\ra P', f\colon M\ra M')$, satisfying the following compatibility conditions
		\begin{enumerate}
			\item  $a_{\mc{G}}'\circ F=a_{\mc{G}}$, and 
			$F(p\cdot \gamma)=F(p)\cdot \gamma$ for all 
			$(p,\gamma)\in P\times_{\mc{G}_0}\mc{G}_1$,
			\item  $\pi'\circ F=f\circ \pi$.
		\end{enumerate}
	\end{definition}
	
	\begin{example}
		Let $\mc{H}=[\mc{H}_1\rra \mc{H}_0]$ be a topological groupoid. Identity assigning map $\mc{H}_0\ra \mc{H}_1$ defines a section of the target map $t\colon  \mc{H}_1\ra \mc{H}_0$. Then $t\colon  \mc{H}_1\ra \mc{H}_0$ can be considered a principal $\mc{H}$-bundle, with action of $\mc{H}$  given by the pair, composition of arrows and the  source map $s\colon  \mc{H}_1\ra \mc{H}_0$.	
	\end{example}
	Given a continuous map $f\colon M'\to M$ between topological spaces, and a principal $\mc{G}$-bundle $(P, \pi, M)$ over $M$, we can obtain a principal $\mc{G}$-bundle $(f^*P, M')$ by the usual pull-back of $(P, \pi, M)$. The action of $\mc{G}$ on $f^*P$ is given by $((p,n),\gamma))\mapsto (\mu(p,\gamma),n)$.
	One of the advantages of working with topological spaces (instead of, let us say, smooth spaces) is that the pullback exists in the category of topological spaces.  Let $\text{Top}$ be the category of (small) topological spaces.  Some authors work with the Cartesian closed category of compactly generated Hausdorff spaces, and they often use the notation $\text{Top}$ for the same.  However, we do not have any such requirement in this paper until Section~\ref{section:Serre-hurewiczfibrations}, and we stick to the conventional meaning of  $\text{Top}$.	In Section~\ref{section:Serre-hurewiczfibrations} we will restrict to the category of compactly generated Hausdorff spaces; the same will be denoted as ${\rm CGTop}$.	
	
	According to our convenience and needs,	we will be using either of the two equivalent definitions of a stack, which we introduce in the following passages. One requires the notion of a fibered category, and the other the notion of a  pseudo-functor. We refer to Chapter $3$ of \cite{MR2223406} to see a detailed correspondence between pseudo-functors and fibered categories. Here we will very briefly describe the correspondence. 
	First, we recall the definition in terms of a pseudo-functor.
	
	\begin{definition}[{\cite{MR2223406}}] Let $\mc{C}$ be a category.
		A \textit{pseudo-functor on $\mc{C}$}, denoted by $\mc{F}\colon \mc{C}^{op}\ra \text{Cat}$, consists of the following data:
		\begin{enumerate}
			\item a category $\mc{F}(U)$ for each object $U$ of $\mc{C}$,
			\item a functor $f^*\colon \mc{F}(V)\ra \mc{F}(U)$ for each morphism $f\colon U\ra V$ of $\mc{C}$,
			\item an isomorphism $\epsilon_U\colon (1_U)^*\Rightarrow 1_{\mc{F}(U)}\colon \mc{F}(U)\ra \mc{F}(U)$ for each object $U$ of $\mc{C}$,
			\item an isomorphism $\alpha_{f,g}\colon (g\circ f)^*\Rightarrow f^*g^*\colon \mc{F}(W)\ra \mc{F}(U)$ for each pair of morphisms $U\xra{f}V\xra{g}W$ of $\mc{C}$,
		\end{enumerate}
		satisfying the following conditions:
		\begin{enumerate}
			\item for every morphism $f\colon U\ra V$ in $\mc{C}$ and every object $\eta$ of $\mc{F}(V)$ we have 
			\[\alpha_{f,1_V}(\eta)=f^*(\epsilon_V\eta), \text{~and~} \alpha_{1_U,f}(\eta)=\epsilon_U(f^*\eta),\]
			\item for every triple of morphisms $U\xra{f}V\xra{g}W\xra{h}T$ in $\mc{C}$ and an object $\eta$ of $\mc{F}(T)$, we have 
			\[\alpha_{gf,h}(\eta)\circ \alpha_{f,g}(h^*\eta)=
			\alpha_{f,hg}(\eta)\circ f^*\alpha_{g,h}(\eta).\]
		\end{enumerate}
	\end{definition}

	Let $B\mc{G}(M)$ denote the category of $\mc{G}$-bundles over the topological space $M$, with arrows being morphisms of $\mc{G}$-bundles of the form $(F, {\rm Id}_M)$.  For a principal $\mc{G}$-bundle $\pi:P\ra M$, the map $\pi\colon P\to M$ admits local sections. That means any morphism of principal $\mc{G}$-bundles over a topological space $M$ of the form $(F,1_M)$ is an isomorphism. In other words, for each $M$ in $\text{Top}$, the category $B\mc{G}(M)$ is an object of $\rm Gpd$, where $\rm Gpd$ be the category of groupoids. 	
	
	Another way to interpret  this  is as a pseudo-functor $B\mc{G}\colon \text{Top}^{\rm op}{\ra \text{Gpd}}\subset \text{Cat}$, which sends an object $M$ in $\rm {Top}$ to the category of $\mc{G}$-bundles over the topological space $M$. Whereas a continuous map $f\colon M'\to M$ is sent to a functor $B\mc{G}(f)\colon B\mc{G}(M)\to B\mc{G}(M')$ defined by the pull-back of principal bundles as follows. Given a principal $\mc{G}$-bundle $P$ over $M$, we obtain the pull-back principal $\mc{G}$-bundle over $M'$. On the other hand,  given an arrow $(F\colon P_1\to P_2, {\rm Id}_M)$ in $B\mc{G}(M)$, the arrow	
	$(f^*P_1\to f^*P_2, {\rm Id}_{M'})$ in $B\mc{G}(M')$ is the unique arrow defined by the universal property of the pull-back diagram with respect to the maps $f\colon M'\to M$  and $P_2\to M$.

	Now we give a description of a stack in terms of a pseudo-functor. 
	
	
	By a \textit{site} we mean a category $\mc{C}$ with a choice of a Grothendieck topology $\mc{J}$ on $\mc{C}$. Let $(\mc{C},\mc{J})$ be a site and $\Phi\colon \mc{C}^{\rm op}\ra \text{Cat}$ be a pseudo-functor.  Let 	$\{\sigma_\alpha\colon U_\alpha\ra U\}$ be a covering of an object $U$ in $\mc{C}$. We simplify our notations as follows. 
	Henceforth $U_{\alpha\beta}$ will denote the fiber product of morphisms $U_\alpha\ra U$ and $U_\beta\ra U$  in the following diagram,
	\begin{equation}\label{Diagram:2pullback}
	\begin{tikzcd}[sep=small]
	U_{\alpha\beta}=U_\alpha\times_UU_\beta \arrow[dd, "{\rm pr}_1"] \arrow[rr, "{\rm pr}_2"] &  & U_\beta \arrow[dd] \\
	&  &                    \\
	U_\alpha \arrow[rr]                                                           &  & U                 
	\end{tikzcd}.
	\end{equation}
	Similarly $U_{\alpha\beta\gamma}$ denote the fiber product of morphisms $U_\alpha\ra U,U_\beta\ra U, U_\gamma\ra U$ as in the following diagram,
	\begin{equation}\label{Diagram:3pullback}
	\begin{tikzcd}[sep=small]
	& U_{\alpha\beta\gamma} \arrow[dd, "{\rm pr}_{13}"] \arrow[rr, "{\rm pr}_{23}"] \arrow[ld, "{\rm pr}_{12}"'] &                    & U_{\beta\gamma} \arrow[dd] \arrow[ld] \\
	U_{\alpha\beta} \arrow[dd] \arrow[rr] &                                                                                          & U_\beta \arrow[dd] &                                       \\
	& U_{\alpha\gamma} \arrow[rr] \arrow[ld]                                                   &                    & U_\gamma \arrow[ld]                   \\
	U_\alpha \arrow[rr]                   &                                                                                          & U                  &                                      
	\end{tikzcd}
	\end{equation}
	For an arrow $f$ in $\mc{C}$ the functor  $\Phi(f)$  will be denoted as $f^{*}$.	
	
	First we associate a category $\Phi_{\text{desc}}(\{U_\alpha\ra U\})$, namely the \textit{descent category} of $\Phi$ with respect to the covering $\{U_\alpha\ra U\}$, as follows. 
	
	\begin{itemize}
		\item An object in $\Phi_{\text{desc}}(\{U_\alpha\ra U\})$ is  a collection  $\big\{\{s_\alpha\},\{\phi_{\alpha\beta}\}\big\}$, where
		\begin{enumerate}
			\item $s_\alpha$ is an object of $\Phi(U_\alpha)$ for each $\alpha\in \Lambda$,
			\item $\phi_{\alpha\beta}\colon {\rm pr}_2^*s_\beta\ra {\rm pr}_1^*s_\alpha$ is an isomorphism in $\Phi(U_{\alpha\beta})$ for each $\alpha,\beta\in \Lambda$,
		\end{enumerate}
		satisfying the condition \[{\rm pr}_{13}^*(\phi_{\alpha\gamma})={\rm pr}_{12}^*(\phi_{\alpha \beta})\circ {\rm pr}_{23}^*(\phi_{\beta \gamma})\]	 in the category $\Phi(U_{\alpha\beta\gamma})$ for each $\alpha,\beta,\gamma\in \Lambda$.
		\item A morphism from $\big\{\{s_\alpha\},\{\phi_{\alpha\beta}\}\big\}$ to $\big\{\{t_\alpha\},\{\psi_{\alpha\beta}\}\big\}$ is given by a collection $\{\theta_\alpha\colon s_\alpha\ra t_\alpha\}_{\alpha\in \Lambda}$. Here each $\theta_\alpha$ is a morphism in $\Phi(U_\alpha)$   for each $\alpha\in \Lambda$ satisfying the condition ${\rm pr}_1^*(\theta_\alpha)\circ \phi_{\alpha\beta}=\psi_{\alpha\beta}\circ {\rm pr}_2^*(\theta_\beta)$  in category $\Phi(U_{\alpha\beta})$ for each $\alpha,\beta\in \Lambda$.
	\end{itemize}
	
	For each object $U$ of $\mc{C}$ and a cover $\{\sigma_\alpha\colon U_\alpha\ra U\}$ of $U$, we have the functor $\Phi(U)\ra \Phi_{\text{desc}}(\{U_\alpha\ra U\})$. The functor is specified for an object $s$ in $\Phi(U)$  as $\{\sigma_{\alpha}^*s, \phi_{\alpha \beta}\}$, where $\phi_{\alpha \beta}$ are isomorphisms derived from the fact that both ${\rm pr}_2^*\sigma_\beta^* s$ and ${\rm pr}_1^*\sigma_\alpha^* s$ are pull-backs of $s$ to $U_{\alpha \beta}$. On the morphism level it is given by $\sigma^*f\colon \sigma^*s\to \sigma^*t$, where $f\colon s\to t$ is an arrow in $\Phi(U)$.
	
	\begin{definition}[\cite{MR2223406}]\label{Def:Stack1}
		Let $(\mc{C},\mc{J})$ be a site. A pseudo-functor $\Phi\colon \mc{C}^{op}\ra \text{Gpd}$ is called a \textit{stack (of groupoids)} over the site $\mc{C}$, if for each object $U$ of $\mc{C}$ and a covering $\{U_\alpha\ra U\}$, the associated functor $\Phi(U)\ra \Phi_{\text{desc}}(\{U_\alpha\ra U\})$ described above is an equivalence of categories.
	\end{definition}
	A pseudo-functor possessing the above property is often called a ``locally determined" pseudo-functor.
	
	Now we turn the category ${\rm Top}$ into a site by specifying a Grothendieck topology and subsequently consider the stack defined over such a site. For our purpose, we will be considering only a very specific type of Grothendieck topology on $\rm {Top}$. 
	Consider the  Grothendieck topology $\mc{J}$ on ${\rm Top}$ described as follows.  A covering of an object $U$ of $\text{Top}$ is given by a jointly surjective family $\{\sigma_\alpha\colon U_\alpha\ra U\}$, where each $U_\alpha \ra U$ is an open immersion. Here, by jointly surjective family, we mean $\bigcup_\alpha \sigma_\alpha(U_\alpha)=U$. We call this topology to be \textit{the open-cover topology on $\text{Top}$}.  Unless mentioned otherwise,  we will always assume $\rm {Top}$ to be equipped with the open-cover topology. 
	
	\begin{example}
		Let $\mc{G}$ be a topological groupoid. Consider the pseudo-functor $B\mc{G}\colon  \text{Top}\ra \text{Gpd}$ mentioned before. Then, 
		$B\mc{G}\colon  \text{Top}\ra \text{Gpd}$ is a stack over the site $(\text{Top},\mc{J})$.
	\end{example}
	
	An alternate way to define a stack is in terms of fibered categories.  In fact, this alternate definition sometimes can be easier to work with.
	For our future purpose, we introduce the definition of a fibered category here.    
	
	Let $\pi_{\mc{D}}\colon \mc{D}\ra \mc{C}$ be a functor. For an object $U$ of $\mc{C}$, let $\mc{D}(U)$ be the category with following description:
	\begin{itemize}
		\item  $\text{Obj}(\mc{D}(U))=\{\eta\in \text{Obj}(\mc{D})|\pi_{\mc{D}}(\eta)=U\}$,
		\item $\text{Mor}_{\mc{D}(U)}(\eta,\eta')=\{f\colon \eta\ra \eta'\in \mc{D}|\pi_{\mc{D}}(f)=1_U\}$.
	\end{itemize}
	We  call  the category $\mathcal{D}(U)$ to be $\textit{the fiber of $U$ in $\mc{D}$}$.			
	
	\begin{definition}[{\cite{MR2223406}}]Let $\pi_{\mc{D}}\colon \mc{D}\ra \mc{C}$ be a functor. A morphism $\theta\colon \xi\ra \eta$ in $\mc{D}$ is said to be a \textit{Cartesian arrow} if, for every morphism $\psi\colon \zeta\ra \eta$ in $\mc{D}$ and a morphism $h\colon \pi_{\mc{D}}(\zeta)\ra \pi_{\mc{D}}(\xi)$ with $\pi_{\mc{D}}(\theta)\circ h=\pi_{\mc{D}}(\psi)$, there exists a unique morphism $\Phi\colon \zeta\ra \xi$ such that $\psi=\Phi\circ \theta$ and $\pi_{\mc{D}}(\Phi)=h$.
		
		A functor $\pi_{\mc{D}}\colon \mc{D}\ra \mc{C}$ is called a \textit{fibered category over $\mc{C}$},	 if for every object $\eta$ of $\mc{D}$ and a morphism $f\colon U\ra V$ in $\mc{C}$, with $\pi_{\mc{D}}(\eta)=f$, there exists a Cartesian arrow $\theta\colon \xi\ra \eta$ in $\mc{D}$ such that $\pi_{\mc{D}}(\theta)=f$.  
		
		A fibered category $\pi_{\mc{D}}\colon \mc{D}\ra \mc{C}$ is called \textit{fibered in groupoids}; if for each object $U$ of $\mc{C}$, the fiber $\mc{D}(U)$ is a groupoid.	
		
		A morphism of fibered categories from $(\mc{D},\pi_{\mc{D}},\mc{C})$ to $(\mc{D}',\pi_{\mc{D}},\mc{C})$ is given by a functor $F\colon \mc{D}\ra \mc{D}'$, which maps a Cartesian arrow  to a Cartesian arrow and satisfies  the condition $\pi_{\mc{D}'}\circ F=\pi_{\mc{D}}$.	
		
	\end{definition}

	
	Let $(\mc{C}, \mc{J})$ be a site and $\pi_{\mc{D}}\colon \mc{D}\ra \mc{C}$ be a category fibered in groupoids. Then we have a pseudo-functor (which we also denote by $\mc{D}!$) 
	$\mc{D}\colon \mc{C}^{\rm op}\to {\rm Gpd}$, which sends an object $U$ of $\mc{C}$ to the groupoid $\mc{D}(U)$, the fiber of $U$ in $\mc{D}$.

	\begin{definition}[\cite{MR2223406}]\label{Def:Stack2}
		Let $(\mc{C}, \mc{J})$ be a site and $\pi_{\mc{D}}\colon \mc{D}\ra \mc{C}$ be a fibered category over $\mc{C}$. 	We say that $\textit{$\mc{D}\ra \mc{C}$ is a stack over the site $(\mc{C},\mc{J})$}$ if the pseudo-functor $\mc{D}\colon \mc{C}^{\rm op}\to {\rm Gpd}$ is a stack (of groupoids) over the site $\mc{C}$.
		
	\end{definition}
	Conversely, given a pseudo functor $\Phi\colon \mc{C}^{\rm op}\ra \text{Gpd}$ as in Definition~\ref{Def:Stack1}, one associates  a category $\mc{D}$ along with a functor $\pi_{\mc{D}}\colon {\mc{D}}\to \mc{C}$ by taking disjoint unions of image of the functor $\Phi$, such that it is a fibered category over $\mc{C}$

	Let $(\mc{S},\mc{J})$ be a site. If the hom functor (also called the functor of points) 
	$h_U\colon \mc{S}^{\rm op}\ra \text{Set}$ is a sheaf on the site $(\mc{S},\mc{J})$ for each object $U$ in 
	$\mc{S}$, then we call $(\mc{S},\mc{J})$ a \textit{subcanonical site} \cite{MR2223406}. 
	\begin{example}\label{Example:everyobjectisstack}
		Suppose $(\mc{S},\mc{J})$  is a subcanonical site. If we treat a set as a category with only unit arrows, then ${\rm Set}$ can be identified as  a subcategory of ${\rm Cat}$ and   $h_U\colon  \mc{S}^{\rm op}\ra \text{Set}\subset \text{Cat}$ defines a pseudo functor for each object $U$ of $\mc{S}$. Moreover this pseudo functor defines a stack on the site $(\mc{S},\mc{J})$. Let $(\mc{S}/U)\ra \mc{S}$ be the corresponding fibered category over the site $(\mc{S},\mc{J})$.	We often denote the stack $(\mc{S}/U)\ra \mc{S}$ by $\underline{U}$ for an object $U$.
	\end{example}

	\begin{definition}[$2$-fiber product of categories fibered in groupoids { \cite{MR2778793}}]\label{Definition:2-fiberproductinCFGs}
		Let $\pi_{\mc{X}}\colon \mc{X}\rightarrow \mc{S}$, $\pi_{\mc{Y}}\colon{\mc{Y}}\rightarrow \mc{S}$ and $\pi_{\mc{Z}}\colon  {\mc{Z}}\rightarrow \mc{S}$ be categories fibered in groupoids. 
		Let    $f\colon \mc{Y}\rightarrow \mc{X}, g\colon \mc{Z}\rightarrow \mc{X}$ be a pair of morphisms of categories fibered in groupoids. We define the \textit{2-fiber product of $ \mc{Y}$ and $ \mc{Z} $ with respect to morphisms $f,g$} to be the groupoid $\mc{Y}\times_{\mc{X}}\mc{Z}$, with 
		\begin{align*}
		(\mc{Y}\times_{\mc{X}}\mc{Z})_0=\big\{(y,z,\alpha)\in \mc{Y}_0\times \mc{Z}_0 \times \mc{X}_1 \mid~ \pi_Y(y)=\pi_Z(z), \alpha\colon f(y)\rightarrow g(z)\big\},\end{align*}
		\begin{align*}		\text{Hom}_{\mc{Y}\times_{\mc{X}}\mc{Z}}
		\big((y,z,\alpha), (y',z',\alpha')\big)&\\
		=\{u\colon y\ra y'\in \mc{Y}_1, v\colon z\ra &z'\in \mc{Z}_1| {\pi_{\mc{Y}}(u)=\pi_{\mc{Z}}(v)}, \alpha'\circ f(u)=g(v)\circ \alpha\}.\end{align*}
		The functor $\pi_{f,g}\colon \mc{Y}\times_{\mc{X}}\mc{Z}\ra \mc{S}$ given by composition $\mc{Y}\times_{\mc{X}}\mc{Z}\xra{pr_1}\mc{Y}
		\xra{\pi_{\mc{Y}}}\mc{S}$ or 
		$\mc{Y}\times_{\mc{X}}\mc{Z}\xra{pr_2}\mc{Z}
		\xra{\pi_{\mc{Z}}}\mc{S}$ turns $\mc{Y}\times_{\mc{X}}\mc{Z}$ into a category fibered in groupoids over $\mc{S}$. We call this category $\mc{Y}\times_{\mc{X}}\mc{Z}$ \textit{the $2$-fiber product of $\mc{Y}$ and $\mc{Z}$ with respect to the morphisms
			$f\colon \mc{Y}\rightarrow \mc{X}$ and
			$g\colon \mc{Z}\rightarrow \mc{X}$}.   \end{definition}
	The definition of $2$-fibered product for categories fibered in groupoids naturally extends to stacks.
	If  $\pi_{\mc{D}}\colon \mc{D}\ra \mc{S}, \pi_{\mc{D}'}\colon \mc{D}'\ra \mc{S}$ and $\pi_{\mc{C}}\colon \mc{C}\ra \mc{S}$ are stacks and, $F\colon \mc{D}\ra \mc{C}$ and $G\colon \mc{D}'\ra \mc{C}$ are morphisms of stacks, then,  $\mc{D}\times_{\mc{C}} \mc{D}'\ra \mc{S}$ is a stack. We call  $\mc{D}\times_{\mc{C}} \mc{D}'\ra \mc{S}$  \textit{the $2$-fiber product stack} of $\mc{D}$ and $\mc{D'}$ with respect to the morphisms $F\colon \mc{D}\ra \mc{C}$ and $F'\colon \mc{D}'\ra \mc{C}$. Note that, for any arbitrary morphisms of stacks $\mc{E}\ra \mc{D}$ and $\mc{E}'\ra \mc{D}$, the stacks $\mc{E}\times_{\mc{D}}\mc{E}'$ and $\mc{E}'\times_{\mc{D}}\mc{E}$ are isomorphic.
	
	\begin{definition}[topological stack 
		{\cite{Carchedi}}]
		\label{Definition:topologicalstack}
		Let $(\text{Top},\mc{J})$ be the site  of topological spaces with open cover topology. A stack $\mc{D}\ra \text{Top}$ is called a \textit{topological stack} if there exists a topological groupoid $\mc{G}=[\mc{G}_1\rra \mc{G}_0]$ and an isomorphism of topological stacks $B\mc{G}\cong \mc{D}$.
	\end{definition}
	\begin{remark}
		In \cite[Definition $7.1$]{noohi2005foundations} the topological stack defined above has been called a pretopological stack. Whereas the topological stack in the same paper has been defined [Definition $13.8$.] as a pretopological stack with some additional conditions.
	\end{remark}
	\begin{example}
		Let $X, Y$ be topological space, and $f\colon   X\ra Y$ be a continuous map. Let $\underline{X}$ and $\underline{Y}$ be stacks associated to $X$ and $Y$ respectively (see \ref{Example:everyobjectisstack}). Then, the map $f\colon  X\ra Y$ induces a morphism of stacks 
		$F\colon  \underline{X}\ra \underline{Y}$  defined by compositions (both at the level of objects and at the level of morphisms)
	\end{example}
	
	Another way we can characterize a topological stack $\mc{D}\ra \text{Top}$ is in terms of an atlas. For this, we need to introduce the following definitions of representable morphisms and epimorphisms of topological stacks. 	
	
	\begin{definition}[representable stack]
		Let $(\mc{S},\mc{J})$ be a subcanonical site. A stack $\mc{D}\ra \mc{S}$ over the site $(\mc{S},\mc{J})$ is said to be a \textit{representable stack} if there exists an object $U$ of $\mc{S}$ and an isomorphism of stacks $\mc{D}\cong \underline{U}$. 
	\end{definition}
	
	
	
	\begin{definition}[representable morphism of stacks {\cite{Carchedi}}]
		Let $(\mc{S},\mc{J})$ be a subcanonical site. Let $\mc{D}$ and $\mc{C}$ be stacks over the site $(\mc{S},\mc{J})$. We call a morphism of stacks  $F\colon  \mc{D}\ra \mc{C}$ a \textit{representable morphism},  if for each object $U$ of $\mc{S}$ and a morphism of stacks $q\colon  \underline{U}\ra \mc{C}$, the $2$-fiber product  stack $\mc{D}\times_{\mc{C}}\underline{U}$ 
		is representable by an object of $\mc{S}$.
	\end{definition}
	\begin{definition}[an epimorphism of stacks {\cite{noohi2005foundations}}]\label{Def:Epistacks}
		Let $(\mc{S},\mc{J})$ be a subcanonical site.
		A morphism of stacks  $F\colon  \mc{D}\ra \mc{C}$ is said to be \textit{an epimorphism of stacks}, if for each object $U$ of $\mc{S}$ and a morphism of stacks $q\colon  \underline{U}\ra \mc{C}$, there exists a covering $\{\sigma_\alpha\colon  U_\alpha\ra U\}$ of $U$ and a family of morphisms $\{q_\alpha\colon  \underline{U_\alpha}\ra \mc{D}\}$ such that the following diagram is $2$-commutative
		\[
		\begin{tikzcd}[sep=small]
		\underline{U_{\alpha}} \arrow[dd, "q_{\alpha}"'] \arrow[rr, "\sigma_{\alpha}"] &  & \underline{U} \arrow[dd, "q"] \\
		&  &                               \\
		\mathcal{D} \arrow[rr,"F"] \arrow[Rightarrow, shorten >=10pt, shorten <=10pt, uurr]                                                   &  & \mathcal{C}                  
		\end{tikzcd},\]
		for each $\alpha$.
	\end{definition}
	
	\begin{definition}
		Let $\mc{D}\ra \text{Top}$ be a stack. Let $X$ be a topological space. A morphism of stacks $f\colon \underline{X}\rightarrow \mc{D}$ is said to be an atlas for the stack $\mc{D}$ if, for every topological space $Y$ and a morphism of stacks $\underline{Y}\ra \mc{D}$, the $2$-fibered product $\underline{X}\times_{\mc{D}}\underline{Y}$ is representable by a topological space $X\times_{\mc{D}}Y$ 
		and the map of topological spaces $X\times_{\mc{D}}Y\rightarrow Y$ associated to the morphism of stacks ${\rm pr}_1\colon \underline{X}\times_{\mc{D}}\underline{Y}\rightarrow \underline{Y}$  is  a map of local sections.
	\end{definition}
	
	An epimorphism of stacks can be described in terms of local sections, which we discuss below. First, let us recall the $2$-Yoneda lemma.
	\subsection{$2$-Yoneda lemma{ \cite{MR2778793}}} 
	
	Let $X$ be an object in the category $\mc{C}$ $\mc{C}$ and $h_X\colon  \mc{C}^{\text{op}}\ra \text{Set}$ be the associated representable functor.  As before we can view this as a pseudo-functor $h_X\colon  \mc{C}^{\text{op}}\ra \text{Cat}$.  Let $\mc{F}\colon  \mc{C}^{op}\ra \text{Cat}$ be a pseudo-functor. Then, the $2$-Yoneda lemma gives an equivalence of categories $\text{Hom}_{\text{Fun}}(h_X,\mc{F})\ra \mc{F}(X)$. Here and onwards we will not differentiate between a morphism $h_X\ra \mc{F}$ and the corresponding object of {the category $\mc{F}(X)$}.
	
	\begin{remark}
		Let  $\mc{D}$ and  $\mc{C}$ be stacks over a subcanonical site $(\mc{S},\mc{J})$.
		Then it is immediate from Definition \ref{Def:Epistacks} that a morphism of stacks  $F\colon  \mc{D}\ra \mc{C}$ is an epimorphism of stacks,  if and only if  for each object $U$ of $\mc{S}$ and an object $b$ of $\mc{C}(U)$, there exists a cover $\{U_\alpha\ra U\}$ of $U$ and objects $a_\alpha\in \mc{D}(U_\alpha)$ such that $F(a_\alpha)\cong b|_{U_\alpha}$ in $\mc{C}(U_\alpha)$ for each $\alpha$. 
	\end{remark}
	
	As we show in the following lemma,  an epimorphism of topological stacks can be described in terms of maps of local sections. 
	\begin{lemma}
		A morphism of topological stacks $F\colon  \mc{D}\ra \mc{C}$ is an epimorphism of topological stacks if and only if, for each topological space $U$ and a morphism of stacks $\underline{U}\ra \mc{C}$ there exists a topological space $V$, a map of local sections
		$\pi\colon  V\ra U$ and a morphism of stacks $\theta\colon \underline{V}\ra \mc{D}$ such that the following diagram is $2$-commutative,
		\[\begin{tikzcd}[sep=small]
		\underline{V} \arrow[dd, "\theta"'] \arrow[rr, "\pi"] &  & \underline{U} \arrow[dd, "q"] \\
		&  &                               \\
		\mathcal{D} \arrow[rr,"F"] \arrow[Rightarrow, shorten >=10pt, shorten <=10pt, uurr]                                                   &  & \mathcal{C}                  
		\end{tikzcd}{\color{red}.}\]
		\begin{proof}
			Suppose $F\colon  \mc{D}\ra \mc{C}$ is  an epimorphism. Let $U$ be a topological space and $q\colon  \underline{U}\ra \mc{C}$ be a morphism of topological stacks. Then, there exists an open cover $\{U_\alpha\ra U\}$ and a morphism of stacks $q_\alpha\colon  \underline{U_\alpha}\ra \mc{D}$ such that the following diagram is $2$-commutative, \[\begin{tikzcd}[sep=small]
			\underline{U_{\alpha}} \arrow[dd, "q_{\alpha}"'] \arrow[rr, "\sigma_{\alpha}"] &  & \underline{U} \arrow[dd, "q"] \\
			&  &                               \\
			\mathcal{D} \arrow[rr,"F"] \arrow[Rightarrow, shorten >=10pt, shorten <=10pt, uurr]                                                   &  & \mathcal{C}                  
			\end{tikzcd},\]
			for each $\alpha\in \Lambda$.
			
			Let $V$ denote the disjoint union $\bigsqcup_{\alpha\in \Lambda} U_\alpha$ and $\pi\colon  V\ra U$ be the map defined as $\pi(x)=x$, if $x\in U_\alpha$ for some $\alpha\in \Lambda$. This map $\pi\colon  V\ra U$ is a continuous map with local sections; there exists a cover $\{U_\alpha\}$ of $U$ and sections $U_\alpha\ra V$ of $\pi\colon  V\ra U$. 
			
			Observe that, each morphism of stacks $\underline{U_\alpha}\ra \mc{D}$ assigns an object of $\mc{D}(U_\alpha)$ by $2$-Yoneda Lemma. So, for each $\alpha\in \Lambda$, we have an object $a_\alpha\in \mc{D}(U_\alpha)$. As $U_\alpha\cap U_\beta=\emptyset$ in $V$, the collection $\{a_\alpha\}_{\alpha\in \Lambda}$ trivially agree on the intersection. Now, as $\mc{D}$ is a stack, this compatible collection $\{a_\alpha\in \mc{D}(U_\alpha)\}$ glue together  to produce an object of $\mc{D}(V)$, which by Yoneda $2$-lemma gives a morphism of stacks $\underline{V}\ra \mc{D}$. So, for the morphism of stacks $\underline{U}\ra \mc{C}$, there exists a map of local sections $\pi\colon  V\ra U$ and a morphism of stacks $\theta\colon  \underline{V}\ra \mc{D}$ with the following $2$-commutative diagram,
			\[\begin{tikzcd}[sep=small]
			\underline{V} \arrow[dd, "\theta"'] \arrow[rr, "\pi"] &  & \underline{U} \arrow[dd, "q"] \\
			&  &                               \\
			\mathcal{D} \arrow[rr,"F"] \arrow[Rightarrow, shorten >=10pt, shorten <=10pt, uurr]                                                   &  & \mathcal{C}                  
			\end{tikzcd}.\]
			
			Conversely, assume that for a morphism of stacks $F\colon  \mc{D}\ra \mc{C}$, and a  morphism of stacks $\underline{U}\ra \mc{C}$ for a topological space $U$, there exists a map of local sections $\pi\colon  V\ra U$ and a morphism of stacks $\theta\colon  \underline{V}\ra \mc{D}$ with a $2$-morphism $F\circ \theta\Rightarrow q\circ \pi$. 
			
			As $\pi\colon  V\ra U$ is a map of local sections, there exists an open cover $\{U_\alpha\}$ of $U$ and sections $\sigma_\alpha\colon  U_\alpha\ra V$ of $\pi\colon  V\ra U$. Thus, we have inclusion $i_\alpha=\pi\circ \sigma_\alpha\colon  \underline{U_\alpha}\ra \underline{V}\ra \underline{U}$ and morphism of stacks $q_\alpha=\theta\circ \sigma_\alpha\colon  \underline{U_\alpha}\ra \underline{V}\ra \mc{D}$. The $2$-morphism $F\circ \theta\Rightarrow q\circ \pi$ induce $2$-morphisms $F\circ q_\alpha\Rightarrow q\circ i_\alpha$ giving following $2$-commutative diagram 
			\[	\begin{tikzcd}[sep=small]
			\underline{U_{\alpha}} \arrow[dd, "q_{\alpha}"'] \arrow[rr, "\sigma_{\alpha}"] &  & \underline{U} \arrow[dd, "q"] \\
			&  &                               \\
			\mathcal{D} \arrow[rr,"F"] \arrow[Rightarrow, shorten >=10pt, shorten <=10pt, uurr]                                                   &  & \mathcal{C}                  
			\end{tikzcd}.\]
		\end{proof}
	\end{lemma}
	\begin{proposition}[{\cite{Carchedi}}]
		A stack $\mc{D}$ over $\text{Top}$ with the open cover topology is a topological stack (Definition \ref{Definition:topologicalstack}) if and only if it has an atlas, that is, a representable epimorphism of stacks $\underline{X}\ra \mc{D}$ for an object $X$ of $\text{Top}$.
	\end{proposition}
	The above proposition gives a direct way to compute the topological groupoid corresponding to a topological stack. 
	
	We end this section by making a note of the following properties of a topological stack, which obviously fail in the smooth set-up. 
	
	\begin{lemma}\label{Lemma:allmapsarerepresentable}
		Let $\mc{X}$ be a topological stack. Then, for every topological space $X$, any morphism of stacks $X\ra \mc{X}$ is representable.
		\begin{proof}
			Let $\mc{X}$ be a topological stack.  That means there exists a topological space $M$ and a morphism of stacks $M\ra \mc{X}$ that is a representable epimorphism.
			
			Let $Y$ be a topological space and $g: Y\ra \mc{X}$  a morphism of stacks. We show that $X\times_{\mc{X}}Y$ is a representable stack. 
			
			As $M\ra \mc{X}$ is representable, for the morphism $X\ra \mc{X}$, the fiber product $X\times_{\mc{X}}M$ is representable and likewise $M\times_{\mc{X}}Y$ is also representable. Considering the morphisms $X\times_{\mc{X}}M\ra M$ and $M\times_{\mc{X}}Y\ra M$, the fiber product $(X\times_{\mc{X}}M)\times_{M}(M\times_{\mc{X}}Y)$ is representable by topological spaces (as all three stacks are representable by topological spaces).  Identifying  $(X\times_{\mc{X}}M)\times_{M}(M\times_{\mc{X}}Y)$ with $X\times_{\mc{X}}Y$, we conclude that the fiber product $X\times_{\mc{X}}Y$ is representable by a topological space. 
		\end{proof}
	\end{lemma}	
	
	As a consequence, we have the following prescription to produce an atlas from another. 
	\begin{lemma}\label{Lemma:nonuniquenessofatlas}
		Let $\mc{D}$ be a topological stack and $X\ra \mc{D}$ be an atlas for $\mc{D}$. Then, for any map of local sections $Y\ra X$, the composition $Y\ra X\ra \mc{D}$ is an atlas for $\mc{D}$.
		\begin{proof}
			As any morphism of stacks from a topological space to a topological stack is a representable morphism of stacks, $Y\ra \mc{D}$ is a representable morphism of stacks. 
			
			Let $M$ be a topological space and $M\ra \mc{D}$ be a morphism of topological stacks. As $X\ra \mc{D}$ is an atlas, the $2$-fiber product $\underline{X}\times_{\mc{D}}\underline{M}$ is representable by a topological space and ${\rm pr}_2:\underline{X}\times_{\mc{D}}\underline{M}\ra \underline{M}$ is a map of local sections,  expressed by the $2$-commutative diagram,
			\[
			\begin{tikzcd}[sep=small]
			\underline{X}\times_{\mc{D}}\underline{M} \arrow[dd] \arrow[rr] &  & \underline{M} \arrow[dd] \\
			&  &              \\
			\underline{X} \arrow[rr]            &  & \mc{D}      
			\end{tikzcd}.\]
			Consider the pullback of ${\rm pr}_1:\underline{X}\times_{\mc{D}}\underline{M}\ra \underline{X}$ along $\underline{Y}\ra \underline{X}$ to obtain the following diagram,
			\[
			\begin{tikzcd}[sep=small]
			\underline{Y}\times_{\underline{X}}\underline{X} \times_{\mc{D}}\underline{M}\arrow[rr] \arrow[dd] &  & \underline{X}\times_{\mc{D}}\underline{M} \arrow[dd] \\
			&  &              \\
			\underline{Y} \arrow[rr]                        &  & \underline{X}           
			\end{tikzcd}.\]
			As $\underline{Y},\underline{X},\underline{X}\times_{\mc{D}}\underline{M}$ are representable by topological spaces, the $2$-fiber product $\underline{Y}\times_{\underline{X}}\underline{X}\times_{\mc{D}}\underline{M}$ is representable by a topological space. 
			
			We identify $\underline{Y}\times_{\underline{X}}\underline{X}\times_{\mc{D}}\underline{M}$ with $\underline{Y}\times_{\mc{D}}\underline{M}$. As $Y\ra X$ is a map of local sections, its pullback $\underline{Y}\times_{\mc{D}}\underline{M}\ra \underline{X}\times_{\mc{D}}\underline{M}$ will also be a map of local sections. In turn, we get the composition $\underline{Y}\times_{\mc{D}}\underline{M}\ra \underline{X}\times_{\mc{D}}\underline{M}\ra \underline{M}$ giving a map of local sections.

			Combining the above two diagrams, we have the following diagram,
			\[
			\begin{tikzcd}[sep=small]
			\underline{Y}\times_{\mc{D}}\underline{M} \arrow[rr] \arrow[dd] &  & \underline{M} \arrow[dd] \\
			&  &              \\
			\underline{Y} \arrow[rr]                        &  & \mc{D}      
			\end{tikzcd}.\]
			Thus, for any topological space $M$ and a morphism of topological stacks $\underline{M}\ra \mc{D}$, the $2$-fiber product $\underline{Y}\times_{\mc{D}}\underline{M}$ is representable by a topological space and the projection map $\underline{Y}\times_{\mc{D}}\underline{M}\ra \underline{M}$ is a map of local sections. Hence, $Y\ra \mc{D}$ is an atlas for $\mc{D}$.
		\end{proof}
	\end{lemma}
	
	\section{Morphism of topological groupoids and topological stacks}\label{Section:morphism of top groupoids and extensions}
	\begin{definition}[{\cite{metzler2003}}]
		Let $\mc{G}=[\mc{G}_1\rra \mc{G}_0]$ and $\mc{H}=[\mc{H}_1\rra \mc{H}_0]$ be topological groupoids. A \textit{$\mc{G}-\mc{H}$ bibundle} consists of \begin{itemize}
			\item a topological space $P$,
			\item a left action of $\mc{G}$ on $P$ given by the pair $(a_{\mc{G}}\colon  P\ra \mc{G}_0, \mu_{\mc{G}}\colon  \mc{G}_1\times_{s,\mc{G}_0,a_{\mc{G}}}P\ra P)$,
			\item a right action of $\mc{H}$ on $P$ given by the  pair $(a_{\mc{H}}\colon  P\ra \mc{H}_0, \mu_{\mc{H}}\colon  
			P\times_{a_{\mc{H}},\mc{H}_0,t}\mc{H}_1\ra P)$,
		\end{itemize}
		satisfying the following conditions:
		\begin{enumerate}
			\item the map $a_{\mc{G}}\colon  P\ra \mc{H}_0$ is a principal $\mc{H}$-bundle,
			\item the map $a_{\mc{H}}\colon  P\ra \mc{H}_0$ is $\mc{G}$-invariant, ; that is $a_{\mc{H}}(\gamma\cdot p)=a_{\mc{H}}(p)$ for all $(\gamma,p)\in \mc{G}_1\times_{\mc{G}_0}P$,
			\item the action of $\mc{G}$ on $P$
			is compatible with the action of $\mc{H}$ on $P$, that is 
			$\gamma\cdot(p\cdot\delta)=(\gamma\cdot p)\cdot\delta$ for each $(\gamma,p,\delta)\in \mc{G}_1\times_{s,\mc{G}_0,a_{\mc{G}}}P\times_{a_{\mc{H}},\mc{H}_0,t}\mc{H}_1$.
		\end{enumerate}
		We denote a $\mc{G}-\mc{H}$-bibundle by $P\colon\mc{G}\to{\mc{H}}$.
	\end{definition}
	\begin{example}\label{Example:principalbundleasbibundle}
		Let $\mc{H}=[\mc{H}_1\rra \mc{H}_0]$ be a topological groupoid and $\pi\colon  P\ra M$ be a principal $\mc{H}$-bundle. The topological groupoid $[M\rra M]$  acts on the topological space $P$ by $(\pi\colon  P\ra M, {\rm pr}_2\colon  M\times_{M}P\ra P)$. Then $P$ is a $[M\rra M]-\mc{H}$-bibundle. We call $P$ to be the $[M\rra M]-\mc{H}$-bibundle associated to the principal $\mc{H}$-bundle $\pi\colon  P\ra M$.
	\end{example}
	\begin{example}\label{Example:morphismofgroupoidsasbibundle}
		Let $(\phi_1,\phi_0)\colon  [\mc{G}_1\rra \mc{G}_0]\ra [\mc{H}_1\rra \mc{H}_0]$ be a morphism of topological groupoids. Recall that, the target map $t\colon  \mc{H}_1\ra \mc{H}_0$ can be considered  a principal $\mc{H}$-bundle. Let $\mc{G}_0\times_{\mc{H}_0}\mc{H}_1$ denote the principal $\mc{H}$-bundle given by the pull back of $t\colon  \mc{H}_1\ra \mc{H}_0$ along  $\phi_0\colon  \mc{G}_0\ra \mc{H}_0$. Then  the map $\text{pr}_1\colon  \mc{G}_0\times_{\phi_0,\mc{H}_0,t}\mc{H}_1\ra \mc{G}_0$ and  the map 
		\begin{eqnarray}
		&&\mc{G}_1\times_{s,\mc{G}_0}(\mc{G}_0\times_{\phi_0,\mc{H}_0,t}\mc{H}_1)\ra \mc{G}_0\times_{\phi_0,\mc{H}_0,t}\mc{H}_1,\nonumber\\
		&&(\gamma,(a,\delta))\mapsto (t(\gamma),\phi_1(\gamma)\circ \delta)\nonumber
		\end{eqnarray}
		define a left action of $\mc{G}$ on $\mc{G}_0\times_{\phi_0,\mc{H}_0,t}\mc{H}_1$. Whereas  $s\circ \text{pr}_2\colon  \mc{G}_0\times_{\phi_0,\mc{H}_0,t}\mc{H}_1\ra \mc{H}_0$ and the map  
		\begin{eqnarray}
		&&(\mc{G}_0\times_{\phi_0,\mc{H}_0,t}\mc{H}_1)\times_{\mc{H}_0,t}\mc{H}_1\ra \mc{G}_0\times_{\phi_0,\mc{H}_0,t}\mc{H}_1,\nonumber\\
		&&((a,\delta),\delta')\mapsto (a,\delta \circ \delta')\nonumber
		\end{eqnarray}
		define a right action of $\mc{H}$ on $\mc{G}_0\times_{\phi_0,\mc{H}_0,t}\mc{H}_1$. Then we have a $\mc{G}-\mc{H}$-bibundle  $\mc{G}_0\times_{\phi_0,\mc{H}_0,t}\mc{H}_1$. We call it the $\mc{G}-\mc{H}$-bibundle associated to the morphism of topological groupoids $(\phi_1,\phi_0)\colon  [\mc{G}_1\rra \mc{G}_0]\ra [\mc{H}_1\rra \mc{H}_0]$ and, denote it by  $\mc{G}_0\times_{\phi_0,\mc{H}_0,t}\mc{H}_1\colon  \mc{G}\ra \mc{H}$.
	\end{example}
	\subsection{composition of bibundles}
	Let $\mc{G}=(\mc{G}_1\rra \mc{G}_0),\mc{H}=(\mc{H}_1\rra \mc{H}_0)$ and $\mc{K}=(\mc{K}_1\rra \mc{K}_0)$ be topological groupoids. Let $P\colon  \mc{G}\ra\mc{H}$ and $Q\colon  \mc{H}\ra \mc{K}$ be bibundles of topological groupoids, described by the  following diagram, 
	\[
	\begin{tikzcd}
	\mc{G}_1 \arrow[dd,xshift=0.75ex,"t"]\arrow[dd,xshift=-0.75ex,"s"'] &                                                          & \mc{H}_1 \arrow[dd,xshift=0.75ex,"t"]\arrow[dd,xshift=-0.75ex,"s"'] &                                                          & \mc{K}_1 \arrow[dd,xshift=0.75ex,"t"]\arrow[dd,xshift=-0.75ex,"s"'] \\
	& P \arrow[ld, "a_{\mc{G}}^P"'] \arrow[rd, "a_{\mc{H}}^P"] &                     & Q \arrow[ld, "a_{\mc{H}}^Q"'] \arrow[rd, "a_{\mc{K}}^Q"] &                     \\
	\mc{G}_0            &                                                          & \mc{H}_0            &                                                          & \mc{K}_0           
	\end{tikzcd}.\]
	Let $P\times_{\mc{H}_0}Q$ denote the pullback of   $a_{\mc{H}}^Q\colon  Q\ra \mc{H}_0$ along the map $a_{\mc{H}}^P\colon  P\ra \mc{H}_0$. The topological groupoid 
	$\mc{H}$ acts on $P\times_{\mc{H}_0}Q$  by $\gamma\cdot  (p,q)=(p\cdot \gamma,\gamma^{-1}\cdot q)$.
	Let $(P\times_{\mc{H}_0}Q)/\mc{H}_1$ be the quotient space.    Consider the maps $\widetilde{a_{\mc{G}}}\colon  (P\times_{\mc{H}_0}Q)/\mc{H}_1\ra \mc{G}_0$  and $\widetilde{a_{\mc{K}}}\colon  (P\times_{\mc{H}_0}Q)/\mc{H}_1\ra \mc{K}_0$ respectively given by $[(p,q)]\mapsto a_{\mc{G}}^P(p)$ and  $[(p,q)]\mapsto a_{\mc{K}}^Q(q)$. Then $(P\times_{\mc{H}_0}Q)/\mc{H}_1$ becomes a $\mc{G}-\mc{K}$-bibundle.
	\begin{definition}[composition of bibundles \cite{Carchedi}]
		The $\mc{G}-\mc{K}$-bibundle $(P\times_{\mc{H}_0}Q)/\mc{H}_1\colon  \mc{G}\ra \mc{K}$ obtained above is called the \textit{composition of bibundles $P$ and $Q$}. 
		
		We denote the composition of bibundles $P$ and $Q$  by $P\circ Q$. 
	\end{definition}
	
	\subsection{morphism of topological stack $B\mc{G}\ra B\mc{H}$ associated to a $\mc{G}-\mc{H}$-bibundle}\label{subsection:morphismassocaitedtobibundle} Consider a bibundle  $Q\colon  \mc{G}\ra \mc{H}$. Let  $\pi\colon  P\ra M$ be a principal $\mc{G}$-bundle over the topological space $M$. In Example \ref{Example:principalbundleasbibundle} we have seen this principal $\mc{G}$-bundle over $M$ can be treated as a $[M\rra M]-\mc{G}$-bibundle. The composition of bibundles $P\colon  [M\rra M]\ra \mc{G}$ and $Q\colon  \mc{G}\ra \mc{H}$ then yields a $[M\rra M]-\mc{H}$ bibundle $P\circ Q\colon  [M\rra M]\ra \mc{H}$. Observe that $[M\rra M]-\mc{H}$ bibundle $P\circ Q\colon  [M\rra M]\ra \mc{H}$ is nothing but a principal $\mc{H}$-bundle over a topological space $M$. This produces a morphism of stacks  $BQ\colon B\mc{G}\ra B\mc{H}$. The morphism of stacks $BQ\colon  B\mc{G}\ra B\mc{H}$ will be called the  morphism of stacks associated to the $\mc{G}-\mc{H}$-bibundle $Q\colon  \mc{G}\ra \mc{H}$. 
	
	\subsection{a morphism of topological stacks associated to a morphism of topological groupoids} 
	Let $(\phi_1,\phi_0)\colon  [\mc{G}_1\rra \mc{G}_0]\ra [\mc{H}_1\rra \mc{H}_0]$ be a morphism of topological groupoids. Now, combining Example \ref{Example:morphismofgroupoidsasbibundle} with  construction in section \ref{subsection:morphismassocaitedtobibundle}, we have a   morphism of topological stacks $B\mc{G}\ra B\mc{H}$ associated to   $(\phi_1,\phi_0)\colon \mc{G}\ra \mc{H}$. This morphism  of stacks $B\mc{G}\ra B\mc{H}$ will be called the morphism of topological stacks associated to the morphism of topological groupoids $\mc{G}\ra \mc{H}$.
	
	\begin{remark}
		It would be pertinent to note here that a  morphism of topological stacks $B\mc{G}\ra B\mc{H}$ may not arise from a morphism of topological groupoids. However, a morphism of topological stacks $B\mc{G}\ra B\mc{H}$ is completely determined by a $\mc{G}-\mc{H}$-bibundle. The smooth version of this observation can be found in \cite[Remark $4.18$]{MR2778793}.
	\end{remark}
	In this paper, we will particularly be interested in a   special type of morphism of stacks, known as a gerbe.
	\begin{definition}[a gerbe over a topological stack {\cite{noohi2005foundations}}] Let $\mc{C}$ be a topological stack. A morphism of topological stacks $F\colon  \mc{D}\ra \mc{C}$  is said to be \textit{a gerbe over the topological stack $\mc{C}$}, if 
		\begin{itemize}
			\item the morphism $F\colon  \mc{D}\ra \mc{C}$ is an epimorphism of topological stacks,
			\item the diagonal morphism $\Delta_F\colon  \mc{D}\ra \mc{D}\times_{\mc{C}}\mc{D}$ is an epimorphism of topological stacks.
		\end{itemize}
	\end{definition}
	With this, we conclude our review of standard material. 	
	\section{Topological groupoid extension and the associated morphism of stacks}\label{section:Topogrpdextensionmorstacks}
	The possible relation between gerbes over differentiable stacks and Lie groupoid extensions in the smooth set-up have been alluded to in  \cite{MR2493616}, without much explicit detail. In \cite{MR4124773} we have given an explicit construction for the correspondence between differentiable gerbes and Lie groupoid extension, with some additional conditions imposed on the differentiable gerbe.  	
	However, in the topological set-up, the notion of a groupoid extension has not been explored.
	
	Let $(F,1_M)\colon  [\mc{G}_1\rra M]\ra [\mc{H}_1\rra M]$ be a topological groupoid extension. As we have seen in section \ref{subsection:morphismassocaitedtobibundle}, corresponding to this morphism of topological groupoid  $F$ we have a morphism of stacks $F:B\mc{G}\to B\mc{H}$. Our main observation is, in particular, when $F$  a topological groupoid extension, the corresponding map of stacks enjoys additional properties of, $F\colon B\mc{G}\to B\mc{H}$  and the diagonal functor $\Delta_F\colon  B\mc{G}\ra B\mc{G}\times_{B\mc{H}}B\mc{G}$, both being epimorphisms of stacks. 	In other words, the morphism of stacks $B\mc{G}\to B\mc{H}$ obtained is a gerbe over a topological stack.

	\begin{theorem}\label{Th:topgrpdextensiontogerbes}
		Let $F\colon  [\mc{G}_1\rra M]\ra [\mc{H}_1\rra M]$ be a topological groupoid extension. Then, the associated morphism of topological stacks $F\colon  B\mc{G}\ra B\mc{H}$ is a gerbe over the topological stack $B\mc{H}$.
		\begin{proof}
			{\bf Step1:} We first prove that $F\colon  B\mc{G}\ra B\mc{H}$ is an epimorphism of topological stacks. 
			
			Let $U$ be an object of $\text{Top}$ and $q\colon  \underline{U}\ra B\mc{H}$ be a morphism of topological stacks. Let $\pi:Q\ra U$ be the principal $\mc{H}$-bundle associated to the morphism of stacks $q\colon  \underline{U}\ra B\mc{H}$. 
			
			Since $\pi:Q\ra U$ admits local sections, there exists an open cover $\{U_\alpha\to U\}_{\alpha\in \Lambda}$ of $U$ and a family of maps $\{\sigma_{\alpha}:U_{\alpha}\to M\}$, such that $\pi|_{\pi^{-1}U_\alpha}\colon  \pi^{-1}(U_\alpha)\ra U_\alpha$ is pullback of the trivial principal $\mc{H}$-bundle $t\colon  \mc{H}_1\ra M$. 
			
			Next, we pullback the principal $\mc{G}$-bundle $t\colon  \mc{G}_1\ra M$ along the morphism $\sigma_\alpha\colon   U_\alpha\ra M$ to obtain  a principal $\mc{G}$-bundle over the topological space $U_\alpha$, for each $\alpha$. In turn we obtain  a morphism of topological stacks $q_\alpha\colon  \underline{U_\alpha}\ra B\mc{G}$. So,  given a morphism of stacks $q\colon  \underline{U}\ra B\mc{H}$, we have an open cover $\{U_\alpha\}_{\alpha\in \Lambda}$ and morphism of stacks $q_\alpha\colon  \underline{U_\alpha}\ra B\mc{G}$ forming the $2$-commutative diagram
			\begin{equation}\label{Diagram:UiUBGBH}
			\begin{tikzcd}[sep=small]
			\underline{U_\alpha} \arrow[dd, "q_i"'] \arrow[rr,"\Phi=\text{inclusion}"] & & \underline{U} \arrow[dd, "q"] \\
			& & \\
			B\mc{G} \arrow[Rightarrow, shorten >=10pt, shorten <=10pt, uurr] \arrow[rr, "F"] & & B\mc{H}
			\end{tikzcd}.
			\end{equation}
			The $2$-commutativity of the diagram follows from the observation that the morphism of stacks $F\colon  B\mc{G}\ra B\mc{H}$ maps the pullback of the principal $\mc{G}$-bundle along a map $\eta\colon  U\ra M$ to the pullback of the principal $\mc{H}$-bundle along a map $\eta\colon  U\ra M$ (\cite[Lemma $5.5$.]{MR4124773}).
			
			We now prove that the diagonal morphism $\Delta_F\colon  B\mc{G}\ra B\mc{G}\times_{B\mc{H}}B\mc{G}$ is an epimorphism of stacks. Observe that, the $2$-fiber product $\mc{G}\times_{\mc{H}}\mc{G}$ is a topological groupoid. As the stackification and Yoneda embedding are preserved under $2$-fiber product, we see that $B(\mc{G}\times_{\mc{H}}\mc{G})\cong B\mc{G}\times_{B\mc{H}}B\mc{G}$. It turns out that the $2$-fibered product $\mc{G}\times_{\mc{H}}\mc{G}$ is a transitive topological groupoid and so is Morita equivalent to a topological groupoid $(K\rra *)$ for some topological group $K$ (\cite[Lemma $5.8$, Lemma $5.9$]{MR4124773}). So, the diagonal morphism $\Delta_F\colon  B\mc{G}\ra B\mc{G}\times_{B\mc{H}}B\mc{H}$ is equivalent to morphism of stacks $B\mc{G}\ra B(K\rra *)$. By similar justification as in the case of $B\mc{G}\ra B\mc{H}$, we can see that $B\mc{G}\ra B(K\rra *)$ is an epimorphism of stacks. Thus, the diagonal morphism $\Delta_F\colon  B\mc{G}\ra B\mc{G}\times_{B\mc{H}}B\mc{G}$ is an epimorphism of stacks. In conclusion, the morphism of topological stacks associated to a topological groupoid extension $\mc{G}\ra \mc{H}$ is a gerbe over the topological stack $B\mc{H}$.
		\end{proof}
	\end{theorem}	
	A detailed proof for the smooth version of the following result is given in \cite{MR4124773}. The proof almost works verbatim in the topological set-up as well. Here we give an outline of the proof.
	\begin{theorem}\label{th:gerbetotopgrpdexten}
		Let $\mc{D},\mc{C}$ be topological stacks and $F\colon  \mc{D}\ra \mc{C}$ be a gerbe over the topological stack $\mc{C}$. Further assume that the diagonal morphism $\Delta_F:\mc{D}\ra \mc{D}\times_{\mc{C}}\mc{D}$ is a representable map of local sections. Then, there exists a topological groupoid extension $\mc{G}\ra \mc{H}$ inducing the morphism $\mc{D}\ra \mc{C}$. Further, this topological groupoid extension is unique upto Morita equivalence.
		\begin{proof}
			As $F\colon  \mc{D}\ra \mc{C}$ is an epimorphism of topological stacks, given an object $U$ of $\text{Top}$ and a morphism of topological stacks $\underline{U}\ra \mc{C}$, there exists a map of local sections $V\ra U$ and a morphism of stacks $\underline{V}\ra \underline{U}$ with the following $2$-commutative diagram, 
			\[\begin{tikzcd}[sep=small]
			\underline{V} \arrow[dd, "p"'] \arrow[rr, "\pi"] &  & \underline{U} \arrow[dd, "\tilde{q}"] \\
			&  &                               \\
			\mathcal{D} \arrow[rr,"F"] \arrow[Rightarrow, shorten >=10pt, shorten <=10pt, uurr]                                                   &  & \mathcal{C}                  
			\end{tikzcd}.\]
			Suppose that $\underline{U}\ra \mc{C}$ is an atlas for the stack $\mc{C}$, then, as $V\ra U$ is a map of local sections, the composition $\underline{V}\ra \underline{U}\ra \mc{C}$ is an atlas for the stack $\mc{C}$ (Lemma \ref{Lemma:nonuniquenessofatlas}). Thus, for the morphism of stacks $F\colon  \mc{D}\ra \mc{C}$ there exists an atlas $q\colon  \underline{V}\ra \mc{C}$ for $\mc{C}$ and a morphism of stacks $p\colon \underline{V}\ra\mc{D}$ with the following $2$-commutative diagram  
			\[
			\begin{tikzcd}[sep=small]
			V \arrow[dd, "p"'] \arrow[rrdd, "q"] &  &        \\
			&  &        \\
			\mc{D} \arrow[rr, "F"]               &  & \mc{C}
			\end{tikzcd}.\]
			The condition that $\mc{D}\ra \mc{D}\times_{\mc{C}}\mc{D}$ is an epimorphism implies that the morphism $p\colon\underline{V}\ra \mc{D}$ is an epimorphism of stacks. On the other hand  the morphism $p:\underline{V}\ra \mc{D}$ is an atlas for the stack $\mc{D}$ follows from  the fact that $\mc{D}\ra \mc{D}\times_{\mc{C}}\mc{D}$ is a representable  map of local sections. The atlases $p:\underline{V}\ra \mc{D}$ and $q:\underline{V}\ra \mc{C}$  respectively produce the topological groupoids $(V\times_{\mc{D}}V\rra V)$ and $(V\times_{\mc{C}}V\rra V)$ representing the stacks $\mc{D}$ and $\mc{C}$ respectively. The morphism of stacks $F:\mc{D}\ra \mc{C}$ will then induce a morphism of topological groupoids $(f,1_V)\colon(V\times_{\mc{D}}V\rra V)\ra (V\times_{\mc{C}}V\rra V)$. The condition that $\Delta_F\colon\mc{D}\ra \mc{D}\times_{\mc{C}}\mc{D}$ is a representable morphism of local sections further impose the condition that $f\colon V\times_{\mc{D}}V\ra V\times_{\mc{C}}V$ is a map of local sections. Thus, we get a topological groupoid extension associated to a gerbe over a topological stack $\mc{D}\ra \mc{C}$ satisfying the extra condition that $\Delta_F$ is a representable morphism of local sections. The Morita invariance follows from the corresponding argument in \cite{MR4124773}.
		\end{proof}
	\end{theorem}
	
	\section{Serre, Hurewicz fibrations and gerbes}\label{section:Serre-hurewiczfibrations}
	
	In the preceding section we have seen that a topological groupoid extension $F\colon  [\mc{G}_1\rra M]\ra [\mc{H}_1\rra M]$ associates a gerbe $F\colon B\mc{G}\ra B\mc{H}$ over a topological stack. In this section, we are going to recall the definition of Serre and Hurewicz stacks introduced in \cite{MR3144243} and study certain properties of a  gerbe over such stacks. 
	
	Consistent with the definition in \cite{MR3144243}, here we restrict to the category $\text{CGTop}$ of compactly generated (Hausdorff) topological spaces with the usual open-cover topology. 
	
	\begin{definition}[\cite{MR3144243}]A topological stack $\mc{D}$ is said to be a \textit{Hurewicz stack} (respectively, \textit{Serre stack}) if it admits a presentation $\mb{X}=[R\rra X]$ by a topological groupoid in which the source (hence also the target) map $s\colon  R\ra X$ is locally a Hurewicz fibration (respectively, Serre fibration).
	\end{definition} 
	
	\begin{definition}[\cite{MR3144243}]
		A morphism of topological stacks $F\colon  \mc{D}\ra \mc{C}$ is said to be a Hurewicz (respectively, Serre) morphism if for every object $U$ of $\text{Top}$ and a morphism of topological stacks $U\ra \mc{C}$, the fiber product $\mc{D}\times_{\mc{C}}U$ is a Hurewicz (respectively, Serre) stack. 
	\end{definition}	
	Here, we say a gerbe $F\colon\mc{D}\ra \mc{C}$ is a {\it{Hurewicz gerbe}} (respectively, a {\it{Serre gerbe}}), if the underline map of stacks is Hurewicz (respectively, Serre).
	
	Let $F\colon  (\mc{G}_1\rra M)\ra (\mc{H}_1\rra M)$ be a 	 topological groupoid extension. Suppose $(\mc{G}_1\rra M)$	is a topological groupoid with  the source map $s\colon  \mc{G}_1\ra M$  a locally a Hurewicz fibration (respectively, Serre fibration). Then $B\mc{G}$ is a  Hurewicz stack (respectively, Serre stack). Following result directly follows from Theorem~\ref{Th:topgrpdextensiontogerbes} and { \cite[Lemma 2.4.]{MR3144243}}. 
	\begin{proposition}
		Let $\mc{G}=(\mc{G}_1\rra M)$ be a topological groupoid with  the source map $s\colon  \mc{G}_1\ra M$ a locally a Hurewicz fibration (respectively, Serre fibration) and $F\colon  (\mc{G}_1\rra M)\ra (\mc{H}_1\rra M)$ a topological groupoid extension. Then the gerbe  $F\colon B\mc{G}\ra B{\mc{H}}$ constructed in Theorem~\ref{Th:topgrpdextensiontogerbes} is a Hurewicz gerbe (respectively, Serre gerbe).
	\end{proposition}

	We give an alternate proof of {\cite[Lemma 2.4.]{MR3144243}}.
	
	\begin{lemma} {\cite{MR3144243}}
		Let $f:\mc{X}\ra \mc{Y}$ be a morphism of topological stacks. If $\mc{X}$ is a Serre or Hurewicz stack, then $f:\mc{X}\ra \mc{Y}$ is a Serre or Hurewicz morphism of stacks.
		\begin{proof} Let $\mc{X}$ be a Serre stack, that is, there exists a topological space $X$ and an atlas $X\ra \mc{X}$ such that the source and target maps of the associated topological groupoid $[X\times_{\mc{X}}X\rra X]$ are (locally) Serre fibrations of topological spaces.
			
			Let $Y$ be a topological space and $g: Y\ra \mc{Y}$ a morphism of topological stacks. We show that the fiber product $\mc{X}\times_\mc{Y}Y$ is representable by a topological groupoid $[M\rra N]$ whose source and target maps are (locally) Serre fibrations of topological spaces.
			
			Let $R$ denote the fiber product $X\times_{\mc{X}}X$, and $R\ra \mc{Y}$ be the composition $R\ra \mf{X}\ra \mc{Y}$. Let $X\ra \mc{Y}$ denote the composition $X\ra \mc{X}\ra \mc{Y}$.  Consider the following pullback diagrams,
			\[
			\begin{tikzcd}[sep=small]
			R\times_{\mc{Y}}Y \arrow[dd] \arrow[rr] &                                         & Y \arrow[dd] \\
			& X\times_{\mc{Y}}Y \arrow[dd] \arrow[ru] &              \\
			R \arrow[rr]                            &                                         & \mc{Y}       \\
			& X \arrow[ru]                            &             
			\end{tikzcd}.\]
			As $\mc{Y}$ is a topological stack, the map $R\ra \mc{Y}$ is representable (Lemma \ref{Lemma:allmapsarerepresentable}), thus, the fiber product $R\times_{\mc{Y}}Y$ is representable by topological space. For the same reason, $X\ra \mc{Y}$ is representable, thus, the fiber product $X\times_{\mc{Y}}Y$ is representable by topological space.  
			
			By pulling back  $s:R\ra X$  along the map $X\times_{\mc{Y}}Y\ra X$ we obtain a map $R\times_{\mc{Y}}Y\ra X\times_{\mc{Y}}Y$. As $s:R\ra X$ is a Serre fibration, the pullback  $R\times_{\mc{Y}}Y\ra X\times_{\mc{Y}}Y$ is a (locally) Serre fibration. Similarly, $t:R\ra X$ (a (locally) Serre fibration of spaces) is pulled back to give a morphism $R\times_{\mc{Y}}Y\ra X\times_{\mc{Y}}Y$ (which would be a (locally) Serre fibration of spaces). Similarly, the other structure maps on $[R\rra X]$ is pulled back to produce a topological groupoid $[R\times_{\mc{Y}}Y\rra X\times_{\mc{Y}}Y]$. It turns out that $\mc{X}\times_{\mc{Y}}Y$ is representable by $[R\times_{\mc{Y}}Y\rra X\times_{\mc{Y}}Y]$. Thus, $\mc{X}\times_{\mc{Y}}Y$ is a Serre stack. The same proof goes verbatim if we consider $\mc{X}$ to be a Hurewicz stack instead of Serre stack. 
		\end{proof}
	\end{lemma}
	An immediate outcome of the above lemma is the following. 
	\begin{corollary}
		Any morphism of stacks $B[X\rra X]\ra B\mc{H}$ induced by a morphism of topological groupoids $[X\rra X]\ra \mc{H}$ is a Serre morphism as well as a Hurewicz morphism. 
	\end{corollary}
	
		\section*{Acknowledgements}
	The first named author acknowledges research support from SERB, DST, Government of India grant MTR/2018/000528. 
	
	\bibliography{topologicalstacks}
	\bibliographystyle{plain}
\end{document}